\documentclass[11pt, twoside]{article}
\usepackage{amsmath,amsthm,dsfont}
\usepackage{amsfonts}
\usepackage{amssymb,latexsym}
\usepackage{enumerate}
\usepackage{blindtext}
\usepackage{parskip}
\newtheorem{theorem}{Theorem}[section]

\newtheorem{proposition}[theorem]{Proposition}
\newtheorem{lemma}[theorem]{Lemma}
\newtheorem{definition}[theorem]{Definition}

\newtheorem{corollary}[theorem]{Corollary}

\def\ind{1{\hskip -3 pt}\hbox{\textsc{I}}}
\def\n{\noindent}
\def\n{\noindent}
\def\fr{\frac}
\def\Om{\Omega}
\def\E{\mathcal E}

\newcommand{\Em}{\mathcal{E}_m}
\def\de{\delta}
\def\u{\tilde{u}}

\def\ve{\varepsilon}
\def\va{\varphi}

\def\F{\mathcal F}
\def\fr{\frac}

\def\O{\Omega}

\def\n{\noindent}

\def\va{\varphi}

\def\de{\delta}

\def\al{\alpha}

\def\R{\mathbb{R}}
\numberwithin{equation}{section}
\DeclareMathOperator{\loc}{loc}

\begin{document}
\setlength{\baselineskip}{18truept}
\pagestyle{myheadings}

\markboth{H. T. Anh, N. V. Phu and N.Q. Dieu}{ Weighted energy class of $m$-subharmonic functions }
\title {Weighted energy class of $m$-subharmonic functions }
\author{
Hoang Thieu Anh*, Nguyen Van Phu** and Nguyen Quang Dieu*** 
\\ *Faculty of Basic Sciences, University of Transport and Communications, 3 Cau Giay,\\ Dong Da, Hanoi, Vietnam.
	\\ **Faculty of Natural Sciences, Electric Power University,\\ Hanoi,Vietnam;\\
	***Department of Mathematics, Hanoi National University of Education,\\ Hanoi, Vietnam;	 
	\\E-mail: anhht@utc.edu.vn, phunv@epu.edu.vn and ngquang.dieu@hnue.edu.vn}

\date{}
\maketitle

\renewcommand{\thefootnote}{}

\footnote{2020 \emph{Mathematics Subject Classification}: 32U05, 32W20.}

\footnote{\emph{Key words and phrases}: $m$-subharmonic functions, Complex $m$-Hessian operator, $m$-Hessian type equations, $m$-polar sets, $m$-hyperconvex domain.}

\renewcommand{\thefootnote}{\arabic{footnote}}
\setcounter{footnote}{0}

\begin{abstract}
In this paper, we introduce the class $\mathcal{E}_{m,F}(\Omega)$ and solve complex $m$-Hessian equations in the class $\mathcal{E}_{m,F}(\Omega)$. Afterthat, we study subextension in the class $\mathcal{E}_{m,F}(\Omega)$ with the weighted Hessian measure of subextension unchanged. This is an extensive version of the result in $\cite{PDmmj}$ in the case when the smaller domain is relatively compact inside the large one.
\end{abstract}

\section{Introduction}
Let $\Omega\subset \tilde{\Omega}$ be bounded domains in $\mathbb{C}^{n}$ and $u\in PSH(\Omega).$ A function $\tilde{u}\in PSH(\tilde{\Omega})$ is said to be a subextension of $u$ if $\tilde{u}\leq u$ on $\Omega.$ In \cite{M80}, El. Mir gave an example of a plurisubharmonic function on the unit bidisc such that the restriction to any smaller bidisc admits no subextension to the whole space. The subextension in the class $\mathcal{F}(\Omega),$ where $\Omega$ is a bounded hyperconvex domain in $\mathbb{C}^{n}$ has been studied by Cegrell and Zeriahi in $\cite{CZ03}.$ In details, the authors proved that if $\Omega\Subset \tilde{\Omega}$ are bounded hyperconvex domains in $\mathbb{C}^{n}$ and $u\in\mathcal{F}(\Omega)$ then there exists $\tilde{u}\in\mathcal{F}(\tilde{\Omega})$ such that $\tilde{u}\leq u$ on $\Omega$ and $\int_{\tilde{\Omega}}(dd^{c}\tilde{u})^{n}\leq \int_{\Omega}(dd^{c}u)^{n}.$ For the class $\mathcal{E}_{p}(\Omega),p>0,$ the subextension was studied by Hiep in $\cite{H08}$. In \cite{W06}, for every hyperconvex domain $\Omega$, J. Wiklund constructed a plurisubharmonic function in the class $\mathcal{E}(\Omega)$ which can not be subextended. After that, L. Hed \cite{H10} constructed a function in the class $\mathcal{N}\smallsetminus\mathcal{F}$ that can not subextend. In \cite{CH08}, R. Czy\.{z} and L. Hed studied subextension of plurisubharmonic functions without increasing the total Monge – Amp\`ere mass and in \cite{HH14} L. M. Hai and N. X. Hong studied subextension of plurisubharmonic functions without changing the total Monge – Amp\`ere mesures.\\
In \cite{Bl1} and \cite{SA12} the authors introduced $m$-subharmonic functions and the complex $m$-Hessian operator $H_m(.) = (dd^c.)^m\wedge \beta^{n-m}$ which is well-defined in the class of locally bounded $m$-subharmonic functions. In \cite{Ch12}, L. H. Chinh introduced the Cegrell classes $\mathcal{F}_m(\Omega)$ and $\mathcal{E}_m(\Omega)$ which are not necessarily locally bounded and  the complex $m$-Hessian operator  is well defined in these classes. Similar to the theory of plurisubharmonic functions, in the theory of $m$-subharmonic functions, subextension problems plays a important role. The subextension in the class $\mathcal{F}_m(\Omega),$ where $\Omega$ is a bounded $m$-hyperconvex domain in $\mathbb{C}^{n}$ has been studied by Le Mau Hai và Trieu Van Dung in $\cite{HD20}.$ In Theorem 3.3 in \cite{PDmmj}, Nguyen Van Phu and Nguyen Quang Dieu studied maximal subextension in the class $\mathcal{E}_{m,\chi}(\Omega).$ Note that the above-mentioned result in \cite{PDmmj} only concerns the total weighted Hessian {\it mass} of subextension. Recently, in \cite{Pjmaa}, the author studied subextension of $m$-subharmonic functions with give bounday value and gave an estimation of the weighted Hessian {\it measure }  of subextension.\\
In this paper, inspired by the research going on in this area, we study subextension in the class $\mathcal{E}_{m,F}(\Omega)$. This is an extensive version of the result in \cite{PDmmj} in the case when the smaller domain is relatively compact inside the large one with the weighted Hessian measure of subextension unchanged. In particular, we partially cover the main result in \cite{HD20} in the case when the smaller domain is relatively compact inside the large one. \\
  The paper is organized as follows. Besides the introduction, the paper has other four sections. In Section 2 we refer the reader to \cite{Pjmaa} for a summary of properties of $m$-subharmonic functions which were introduced and investigated intensively in recent years by many authors (see \cite{Bl1}, \cite{SA12}, \cite{Ch15} and \cite{T19}). In the final part of this Section, we recall  the comparison principle for the operator $H_{m,F}.$  In Section 3, we study the convergence of weighted $m$-Hessian measure in $Cap_m$-capacity. In Section 4, firstly we introduce the class $\mathcal{E}_{m,F}(\Omega)$ and research some properties related to the class. Afterthat, we solve complex $m$-Hessian type equations $-F(u,z)H_{m}(u)=\mu$ in the class $\mathcal{E}_{m,F}(\Omega)$ in the case when measure $\mu$ vanishes on all of $m$-polar sets and in the case when  measure $\mu$ is arbitrary. Finally, in Section 5, we study subextension in the class $\mathcal{E}_{m,F}(\Omega).$ \\
\section{Preliminaries}
The elements of the theory of $m$-subharmonic functions and the complex $m$-Hessian operator can be found e.g. in \cite{Bl1}, \cite{Ch12}, \cite{Ch15}, \cite{SA12} and \cite{T19}. A summary of the properties required for this paper can be found in Preliminaries Section (from subsection 2.1 to subsection 2.6) in \cite{Pjmaa}.\\
 \n Now we recall the definition of $\tilde{C}_m(E)$-capacity in \cite{HP17} as follows.
\begin{definition}
	Let $E$ be a Borel subset. We define
$$\tilde{C}_m(E)= \tilde{C}_m(E,\Omega)= \sup\Bigl\{\int\limits_{E}H_{m-1}(u): u\in SH_m(\Omega), -1\leq u\leq 0\Bigl\}.$$
\end{definition}
\n According to the argument after Definition 2.7 in \cite{HP17}, we have $ \tilde{C}_m(E)\leq A(\Omega) {Cap}_m(E).$ It follows that if $u_j$ converges to $u$ in $Cap_m$-capacity then $u_j$  also converges to $u$ in $\tilde{C}_m$-capacity.\\
 \n As in the plurisubharmonic case (see \cite{BT1}), we have the following quasicontinuity property of $m$-subharmonic functions (see Theorem 2.9 in \cite{Ch12} and Theorem 4.1 in \cite{SA12}).
 \begin{proposition} \label{capacity}
 	Let $u \in SH_m (\Om)$. Then for every $\ve>0$ we may find an open set $U$ in $\Om$ with $Cap_m (U)<\ve$
 	and $u|_{\Omega \setminus U}$ is continuous.	
 \end{proposition}
 \n
 Using the above result, as in the plurisubharmonic case (see \cite{BT1}), we have the following important fact about negligible sets for $m$-subharmonic functions (see Theorem 5.3 in \cite{SA12}).
 \begin{proposition} \label{negli}
 	Let $\{u_j\}$ be a sequence of negative $m$-subharmonic functions on $\Om.$ Set 
 	$u:=\sup\limits_{j \ge 1} u_j$. Then the set $\{z \in \Om: u(z)<u^* (z)\}$ is $m-$polar.	
 \end{proposition} 
\begin{definition} {\rm 
Let $\mathcal{K}$ be one of the classes $\mathcal{E}_{m}^{0}(\Omega), \mathcal{F}_{m}(\Omega), \mathcal{N}_{m}(\Omega), \mathcal{E}_{m}(\Omega).$ Denote by $\mathcal{K}^{a}$ being the set of all function in $\mathcal{K}$ whose Hessian measures vanish on all $m-$polar sets of $\Omega$. }
\end{definition}

\begin{definition} {\rm 
 Let $F: \mathbb{R}^{-}\times \Omega\rightarrow \mathbb{R}^{-}$ be such that $F(.,z)$  is a nondecreasing function for all $z\in\Omega$. For each $u\in\mathcal{E}_{m}(\Omega)$, we put 
$$H_{m,F}(u)= -F(u(z),z){(dd^{c}u)}^{m}\wedge \beta^{n-m}.$$}
\end{definition}
\n We recall a version of the comparison principle for the operator $H_{m,F}$ (see Theorem 3.8 in \cite{PDjmaa}) .
\begin{theorem}\label{comparison}
	Suppose that the function $t \mapsto F(t,z)$ is nondecreasing in $t.$ 
	Let $u\in\mathcal{N}_{m}(f), v\in\mathcal{E}_{m}(f)$ be such that $H_{m,F}(u)\leq H_{m,F}(v).$ 
	Assume also that $H_m (u)$ puts no mass on $m$-polar sets.
	Then we have $u\geq v$ on $\Om.$
\end{theorem}

\section{Convergence of weighted complex Hessian operator in $Cap_m$-capacity}
We start with the following lower semi-continuity of measures.
\begin{proposition}\label{weakcap}
	Let $F: [-\infty,0] \times \Om \to [-\infty, 0]$ be an upper semicontinuous function such that for 
	each $z \in \Om$, the function $t \mapsto F(t,z)$ is increasing.
	Let $\{u_j\}_{j \ge 1}, u$ be negative subharmonic functions on $\Om$ and $\{\mu_j\}_{j \ge 1}, \mu$ be Radon measures 
	on $\Om.$ Assume that the following conditions hold true:
			
	\n 
	(i) $u_j\longrightarrow u$ in Lebesgue measures on $\Om;$
	
	\n 
	(ii) $\mu_j \to \mu$ weakly on $\Om.$
	
	Then we have 
	$$ \liminf_{j\to\infty} [-F(u_j,z)\mu_j] \ge  -F(u,z)\mu.$$
\end{proposition}
\begin{proof} 
Fix $f\in C_0^\infty (\O)$ with $0\leq f\leq 1.$ We have to prove that
\begin{equation} \label{eq00}
\liminf_{j\to\infty} \int_{\O} -fF(u_j,z) d\mu_j \ge \int_{\O} -fF(u,z)d\mu.\end{equation}	
Towards this end, for $j \ge 1,$ we set 
$$F_j (t,z):= \sup \{ F(t',z')-j (|t-t'|+\Vert z-z'\Vert^2)^{1/2}: (t',z') \in  \R^- \times \Om\}.$$
It is not hard to check that $F_j \searrow F, F_j$ is continuous, and moreover
$t \mapsto F_j (t,z)$ is increasing for every $j \ge 1$ and $z \in \Om$.
We now show that
	\begin{equation}\label{eq0}\liminf_{j\to\infty} \int_{\O} -fF_j (u_j,z) d\mu_j\geq\int_{\O} -fF (u,z)d\mu.\end{equation}	
Indeed, choose a sequence of functions $\{\varphi_k\}\subset  C({\Omega})$ such that $\varphi_k\searrow u$  on $K:= \text{supp}\ f.$
	By (i) we have $u_{j}\to u$ in $\mathcal{D}^{'}(\Omega).$ 
	So using Theorem 3.2.13 in \cite{Ho94} we get
	$$ \limsup_{j\to\infty}\sup\limits_{K}(u_{j}-\varphi_{k})\leq \sup\limits_{K} (u-\varphi_{k})< 0.$$
	Thus, fix $k\in \mathbb N^*$,  we can find $j_0 \ge k$ such that $u_j\leq \varphi_k$ on $K$ for all $j\geq j_0$. Hence, by the assumption (ii) we get  the following estimates	
	$$\begin{aligned} 
	\liminf_{j\to\infty} \int_{\O} -fF_j (u_j,z) d\mu_j &\geq \liminf_{j\to\infty} \int_{\O} -fF_j (\va_k,z) d\mu_j\\
	&\geq \liminf_{j\to\infty} \int_{\O} -fF_k (\va_k,z) d\mu_j\\
	&=\int\limits_{\O} -fF_k (\va_k,z) d\mu.
	\end{aligned}$$
	Now for $l \ge 1$, by monotone convergence theorem we obtain
	$$\liminf_{k\to\infty} \int\limits_{\O} -fF_k (\va_k,z) d\mu \ge \liminf_{k\to\infty} \int\limits_{\O} -fF_l (\va_k,z) d\mu =\int\limits_{\O} -fF_l (u,z) d\mu.$$
	It follows that 
	$$\liminf_{j\to\infty} \int_{\O} -fF_j (u_j,z) d\mu_j \ge 	\int\limits_{\O} -fF_l (u,z) d\mu.$$
Now by letting $l \to \infty$ and applying again monotone convergence theorem we get (\ref{eq0}). Combining this estimate with the obvious one
$$\liminf_{j\to\infty} \int_{\O} -fF(u_j,z) d\mu_j \ge \liminf_{j\to\infty} \int_{\O} -fF_j (u_j,z) d\mu_j$$
we finish the proof of (\ref{eq00}).\end{proof}

\begin{theorem}\label{cap} Let $F: [-\infty, 0] \times \Om \longrightarrow\mathbb (-\infty,0]$ be a continuous function such that for $z \in \Om,$ the function $t \mapsto F (t,z)$ is increasing in $t.$ 
Let $\{u_j\}_{j \ge 1}, u$ be $m$-subharmonic functions on $\Om$ and $\{\mu_j\}_{j \ge 1}, \mu$ be Radon measures 
	on $\Om.$ Assume that:
	
	\n 
	(i) $u_j\longrightarrow u$ in ${Cap}_{m}$-capacity;
	
	\n 
	(ii) $\int\limits_{P \cap \{u>-\infty\}} d\mu=0$ for all $m$-polar subsets $P$ of $\Om;$
	
	\n 
	(iii)  $F (-\infty, z)h\mu_j \to F (-\infty, z)h\mu$ weakly on $\Om$ for all $h\in SH_m(\Omega)\cap L_{\loc}^{\infty}(\Omega).$
	
	\n
	Then $-F (u_j,z) \mu_j \to -F (u,z)\mu \ \text{weakly}$ on $\Om.$
\end{theorem}
\begin{proof} 
	Fix $f\in C_0^\infty (\O)$ with $0\leq f\leq 1$. 
	By (i) we have  $u_j\to u$ in the Lebesgue measure and by (iii)
	$\mu_j \to \mu$ weakly.
	So using Proposition \ref{weakcap} we get
	\begin{equation}\label{eq1}\liminf_{j\to\infty} \int_{\O} -fF (u_j, z) d\mu_j\geq\int_{\O} -fF (u,z)d\mu.\end{equation}
	Next, we prove that
	\begin{equation}\label{eq2}\limsup_{j\to\infty} \int_{\O} -fF (u_j,z)d\mu_j\leq\int_{\O} -fF(u,z)d\mu.\end{equation}
	By considering $-\fr{F(t,z)}{F(-\infty,z)}$ instead of $F (t,z)$, using (iii), we may achieve that
$F (-\infty,z)=-1$ for all $z \in \Om.$
So in particular
$$F (t,z) \ge -1, \ \forall (t, z) \in [-\infty, 0) \times \Om.$$
	Fix $k\in\mathbb N^*$ and $\ve>0$. Since $F$ is bounded from below and continuous at $-\infty$, we can find $\ve'$ such that $\forall z \in K:= \text{supp}\ f$ and $u(z)>-\infty$ we have
	$$0 \le F (u(z),z)-F (u(z)-\ve',z)<\ve.$$
	Furthermore, by quasicontinuity of $u$ and the fact that $u_j \to u$ in $Cap_m$-capacity,
	we can find an open subset $G_k\supset\{u=-\infty\}$ of $\Omega$ such that
	
	\n 
	(a) $Cap_m(G_k)<\frac 1 {2^k}, u$  is continuous on $\Omega\backslash G_k$;
	
	\n 
	(b) $u_j \ge u-\ve'$ on $K \backslash G_k.$ 
	
	\n 
	Since $G_{k}$ is an open subset of $\Omega$ we have $$h_k:=h_{m,G_{k},\Omega}=h^{*}_{m,G_{k},\Omega} \in SH^{-}_{m}(\Omega).$$ 
	Set $$\tilde G_k:= \bigcup\limits_{j \ge k} G_j, \tilde h_k:= h_{m, \tilde G_k, \Om}.$$ 
	Then we have $h_k \ge \tilde h_k$ and 
	$$\int\limits_{\Om} H_m (\tilde h_k)= C_m(\tilde G_k) \le \sum_{j \ge k} C_m(G_j) \le \sum_{j \ge k} 2^{-j}=2^{1-k}.$$
	So by monotonicity of $F$ and (b) we get
	\begin{align}\label{e3.5}
		\begin{split}
		\int_{\O} -fF(u_j, z)d\mu_j&=\int\limits_{\O\backslash G_k} -fF (u_j,z)d\mu_j+\int_{G_k} -fF(u_j,z)d\mu_j\\
		&\leq\int\limits_{\O\backslash G_k} -fF(u-\ve',z)d\mu_j+\int_{G_k} fd\mu_j\\
		\\&\leq \int\limits_{\O\backslash G_k} -fF (u,z)d\mu_j+\ve \int\limits_{\O} f d\mu_j
		+\int\limits_{\O} -fh_k d\mu_j.
		\end{split}
	\end{align}
\n  Note that since $K$ is a compact subset, we have $\ind_{K}$ is an upper semicontinuous function and since $G_k$ is an open subset, we have $\ind_{G_k}$ is a lower semicontinuous function. Thus, by Lemma 1.9 in \cite{De93} we have
\begin{align*}
	\lim\limits_{j\to\infty}\int\limits_{\O\backslash G_k} -fF (u,z)d\mu_j&=\lim\limits_{j\to\infty} \int\limits_{\O} -fF (u,z)d\mu_j -\lim\limits_{j\to\infty}\int\limits_{ G_k} -fF (u,z)d\mu_j\\
	&= \lim\limits_{j\to\infty}\int\limits_{K} -fF (u,z)d\mu_j -\lim\limits_{j\to\infty}\int\limits_{ G_k} -fF (u,z)d\mu_j\\
	&\leq \int\limits_{K} -fF (u,z)d\mu-\int\limits_{ G_k} -fF (u,z)d\mu\\
	&\leq \int\limits_{\Omega} -fF (u,z)d\mu-\int\limits_{ G_k} -fF (u,z)d\mu\\
	&\leq \int\limits_{\O\backslash G_k} -fF (u,z)d\mu.
\end{align*}
	So using again the assumption (iii) in letting $j\to\infty$ in inequality \eqref{e3.5} we obtain
	$$\begin{aligned} 
	\limsup_{j\to\infty}\int_{\O} -fF (u_j,z)d\mu_j &\leq \int\limits_{\O\backslash G_k} 
	-fF(u,z)d\mu+\ve \mu (K)+\int\limits_{\O} -fh_{k} d\mu\\
	&\leq \int\limits_{\O\backslash \{u=-\infty\}} 
	-fF (u,z)d\mu+\ve \mu(K)+\int\limits_{\O} -f\tilde h_{k} d\mu\\
	\end{aligned}$$
	Here the last inequality follows from the choice of $G_k$ and the fact that $h_k \ge \tilde h_k$ on $\Om.$	\\
	Set $\psi:= (\lim\limits_{k \to \infty} \tilde h_k)^*$. By the comparison principle we obtain
	$$\int_{\Om} H_m (\psi) \le \int\limits_{\Om} H_m (\tilde h_k) \to 0 \ \text{as}\ k \to \infty.$$
	Hence $\psi=0$ on $\Om$ and so $\tilde h_k \uparrow 0$ on $\Om \setminus P$ where $P$ is a $m$-polar set.
	So by letting $k\to\infty$ in the last estimate and making using of (ii) we get
	$$\begin{aligned}
	\limsup_{j\to\infty}\int_{\O} -fF (u_j,z)d\mu_j
	&\leq\int\limits_{\O\backslash\{u=-\infty\}} -fF (u,z)d\mu+\ve \mu (K)+ 
	\int\limits_{P} fd\mu\\
	&\leq\int\limits_{\O\backslash\{u=-\infty\}} -fF (u,z)d\mu+\int\limits_{P \cap \{u=-\infty\}} fd\mu+\\
	&+\int\limits_{P \cap \{u>-\infty\}} fd\mu+\ve \mu (K)
	\\
	& \leq \int_{\O} -fF (u,z)d\mu+\ve \mu (K).
	\end{aligned}$$
	Here the last estimate follows from the facts that $F(-\infty,z)=-1$ for $z \in \Om.$
	By letting $\ve \to 0$ we then get \eqref{eq2}. Putting \eqref{eq1} and  \eqref{eq2}
	together we obtain
	$$\lim_{j\to\infty} \int_{\O} -fF (u_j,z)d\mu_j=\int_{\O} -fF (u,z)d\mu.$$
	Hence, $-F (u_j,z)d\mu_j\to -F (u,z)d\mu$ weakly as we want.
\end{proof}
\n 
As a by product of the above results, we get the following variation of Theorem 3.8 in \cite{HP17}.
\begin{corollary}\label{nicecap}
	Let $F: [-\infty, 0] \times \Om \longrightarrow\mathbb (-\infty,0]$ be a continuous function 
such that $t \mapsto F (t,z)$ is increasing for all $z \in \Om.$
Let $\{u_j, u\}\subset\mathcal E_{m}(\Omega)$, be such that $u_j\geq v$, $\forall j\geq 1$ for some $v\in\mathcal E_{m}(\Omega)$ and that $u_j \to u \in \mathcal E_{m}(\Omega)$ in $Cap_m-$ capacity. Then 
	$-F (u_j,z)H_{m}(u_{j})\to -F (u,z)H_{m}(u)$ weakly.
\end{corollary}	
\begin{proof}
	Set $\mu_j:= H_m (u_j), \mu:= H_m (u).$ By Theorem 3.9 in \cite{HP17} we have 
	$h\mu_j \to h\mu$ weakly on $\Om$ for every locally bounded $m$-subharmonic function $h$. 
Now we claim that $F (-\infty, z)h\mu_j \to F (-\infty, z)h\mu$ weakly on $\Om.$
For this, it suffices to check the convergence on a small ball $B \Subset \Om.$ Since we may approximate $F (-\infty,z)$ uniformly on $B$ by (real) polynomials and since $\sup_{j \ge 1}\mu_j (B)<\infty,$ 
we may suppose $F (-\infty,z) \in \mathcal C^{\infty}_0 (\Om).$
So by Lemma 3.10 in \cite{Ch15}, there exists {\it plurisubharmonic} functions $h_1, h_2\subset \mathcal{E}_m^{0}(\Omega)\cap C(\Om)$ such that $F (-\infty, z)=h_1(z)-h_2(z)$ for all $z \in B.$ By Proposition 3.7 in \cite{HP17}, we infer that $hF(-\infty,z) =h_3-h_4$ where $h_3,h_4\in SH_m(\Omega)\cap L_{\loc}^{\infty}(\Omega).$ Thus applying Theorem 3.9 in \cite{HP17} to $h_3, h_4$ we may conclude that 
$F (-\infty, z)h\mu_j \to F (-\infty, z)h\mu$ weakly on $\Om$ as desired.
Finally, since $H_m(u)$ puts no mass on $P \cap \{u>-\infty\}$ for all $m-$polar sets $P$, by Theorem \ref{cap},
we finish the proof.	
\end{proof}	

\section{Complex $m$-Hessian equations in the class $\mathcal{E}_{m,F}(\Omega)$}
\subsection{The class $\mathcal{E}_{m,F}(\Omega)$}
Throughout this section, we 
let $F: [-\infty, 0] \times \Om \longrightarrow\mathbb [-\infty,0]$ be an upper-semicontinuous function satisfying the following condition:

\n
(A) For $z \in \Om, t \mapsto F (t,z)$ is increasing in $t.$ 

\n Taking the idea from Definition 3.1 in \cite{DT23}, we define
\begin{align*}
	&\mathcal{E}_{m,F}(\Omega):=\{u\in SH_{m}^{-} (\Omega): \ \exists(u_{j})\in\mathcal{E}_{m}^{0}(\Omega), u_{j}\searrow u\,\,\text{such that}\\
	&\hskip4cm \sup_{j}\int_{\Omega}[-F(u_j(z),z)] H_{m}(u_{j})<+\infty\}.
\end{align*}
\n Under further restriction on $F$ we have the following inclusion of weighted energies classes.
\begin{theorem}\label{theo4.1}
 Assume that for all $t<0$ given, there exists a constant $\de_t<0$ (depend only on $t$) such that  \\
\n (B)  $\sup\limits_{z\in\Omega}F(t,z)=\delta_t<0.$ \\
	Then we have $\mathcal{E}_{m,F}(\Omega)\subset \mathcal{N}_m(\Omega).$
\end{theorem}
\begin{proof} By considering $\max \{F, -2\}$ instead of $F$ we may suppose further that 	
	$F \ge -2$ on $[-\infty, 0] \times \Om.$\\
	Let $u\in \mathcal{E}_{m,F}(\Omega)$ and $(u_{j})\in\mathcal{E}_{m}^{0}(\Omega)$ decreasing to $u$ such that 
	$$\sup_{j}\int_{\Omega}[-F(u_j(z),z)] H_{m}(u_{j})<+\infty.$$
 So for all $j$ we have
	$$-F(u_{j}(z),z)H_{m}(u_{j})\leq  H_{m}(2^{1/m} u_{j}).$$
	By Theorem 5.9 in \cite{HP17}, there exists $w_{j}\in\mathcal{E}_{m}^0 (\Omega), w_{j}\geq 2^{1/m} u_{j}$ such that 
	\begin{equation} \label{eq21}
	-F(u_{j}(z),z)H_{m}(u_{j}) = H_{m}(w_{j}).
\end{equation}
Fix $t<0.$ We will show that
	$$u_{j}(z)\geq h_j (z):= \frac{w_{j}(z)}{(-\delta_t)^\frac{1}{m}}+t.$$
Indeed, the estimate is trivial	 on the set $\{u_{j}\geq t\}$ since $0\geq w_{j}$ there.\\
On the other hand, on $\Om_t:= \{u_{j}< t\} \Subset \Om$,  we have 
$$-F (u_j, z) \ge -F(t,z)\geq -\delta_t, \ \forall z \in \Om.$$ 
It follows from (\ref{eq21}) that
$$-\delta_t H_m (u_j) \le H_m (w_j)=-\delta_t H_m (h_j).$$ This is equivalent to $$H_m(u_j)\leq H_m(h_j).$$
We also have $$\liminf\limits_{z\to\partial \Omega_{t}}u_{j}(z) \ge t\geq \limsup\limits_{z\to\partial \Omega_{t}}h_j (z) .$$
Thus, by the comparison principle (Theorem 2.13 in \cite{Ch15}) we obtain $u_{j}(z)\geq h_j(z)$ on 
$\Omega_{t}.$ So we have 
		$$u_{j}(z)\geq  \frac{w_{j}(z)}{(-\delta_t)^\frac{1}{m}}+t.$$
	Now we put $v_j=(\sup\limits_{k\geq j}w_k)^{*}.$ Then we have $v_j\in\mathcal{E}_m^{0}(\Omega)$ and $v_j\searrow :=v$ as $j\to +\infty.$
Since $v_j \ge w_j,$ by the comparision principle we also have
	$$\sup\limits_{j\geq 1}\int\limits_{\Omega}H_{m}(v_{j})
		\leq \sup\limits_{j\geq 1}\int\limits_{\Omega}H_{m}(w_{j})
		=\sup\limits_{j\geq 1}\int\limits_{\Omega}[-F(u_{j}(z),z)]H_{m}(u_{j})<\infty.$$
	Thus, we get $v\in\mathcal{F}_{m}(\Omega).$
	Note that for every $k\geq j$ we have
	$$u_{j}\geq u_{k}\geq \frac{w_{k}(z)}{(-\delta_t)^\frac{1}{m}}+t.$$
	Thus, by taking supremum of right-hand side with respect to $k$,  we obtain 
	$$u_{j}(z)\geq \frac{v_{j}(z)}{(-\delta_t)^\frac{1}{m}}+t.$$
	Let $j\to\infty, $we obtain 
	$$u(z)\geq \frac{v(z)}{(-\delta_t)^\frac{1}{m}}+t.$$
	This implies that  
	$$\tilde{u}(z)\geq \frac{\tilde{v}(z)}{(-\delta_t)^\frac{1}{m}}+t,$$
	where $\tilde{u}$ and $\tilde{v}$ as in Definition of the class $\mathcal{N}_m(\Omega)$. Since $v\in\F_m(\Om),$ we infer that $\tilde{v}=0.$
	Thus, we obtain that $\tilde{u}\geq t.$
	By letting $t\to 0$ we conclude that $\tilde{u} \ge 0$, so $u\in\mathcal{N}_{m}(\Omega)$ as well.
	\end{proof}

\n Next, we have the following simple fact about finiteness of $m$-weighted energies in $\mathcal{E}_{m,F}(\Omega)$.
\begin{proposition}\label{theo4.4}
	Let $F: [-\infty, 0] \times \Om \longrightarrow\mathbb [-\infty,0]$ be an upper-semicontinuous function
	that satisfies the conditions (A) and (B). Assume that $u\in\mathcal{E}_{m,F}(\Omega)$. Then we have $$\int\limits_\Omega[-F(u(z),z) ]H_m(u)<\infty.$$
\end{proposition}
\begin{proof}
	Let ${u_j}$ be a sequence in $\mathcal{E}_m^0(\Omega)$ decreasing to $u$ such that
	$$M:=\sup\limits_j\int\limits_\Omega[-F( u_j(z),z)]H_m(u_j)<\infty.$$
It follows from $u\in\mathcal{E}_{m,F}(\Omega)$ and Theorem \ref{theo4.1} that $u\in \mathcal{N}_m(\Omega).$ By Theorem 3.8 in $\cite{HP17}$, we have $H_m (u_j)\rightharpoonup H_m (u)$. Since $u_j \to u$ in Lebesgue measure, by Proposition \ref{weakcap} we obtain
	$$\int\limits_\Omega[-F( u(z),z)]H_m(u)\le 
	\liminf\limits_{j\rightarrow\infty}	\int\limits_\Omega[-F( u_j(z),z)]H_m(u_j) \le  M<\infty.$$
	The proof is completed.
\end{proof}
\n 
\begin{theorem}\label{theo4.2}
	Let $F: [-\infty, 0] \times \Om \longrightarrow\mathbb [-\infty,0]$ be a continuous function
	that satisfies the conditions (A) and (B). 
	Assume that $u\in\mathcal{N}_m(\Omega)$ and $$\int\limits_\Omega[-F (u(z),z)]H_m(u)<+\infty.$$
	Then there exists a sequence $\{u_j\}\in\E_m^0(\Omega)$ such that $u_j\searrow u$ and 
	$$\lim\limits_{j\rightarrow\infty}\int\limits_\Omega[-F( u_j(z),z)]H_m(u_j)=\int\limits_\Omega[-F( u(z),z)]H_m(u).$$
\end{theorem}

\begin{proof} 
We let $\rho\in\E_m^0(\Omega)\cap C^\infty(\Omega)$ be a defining function for $\Omega$. For each $j$, by Lemma 5.5 in \cite{T19}, it follows that
	$$\mathds{1}_{\{u>j\rho\}}H_m(u)(\Omega)=\int\limits_{\{u>j\rho\}}H_m(u)\leq \int\limits_{\{u>j\rho\}}H_m(j\rho)<\infty.$$
	Obviously, we have
	$$H_m(u)\geq \mathds{1}_{\{u>j\rho\}}(dd^cu)^m\wedge\omega^{n-m}.$$
	By Theorem 5.9 in \cite{HP17}, it follows that there exists $u_j\in\Em(\Omega)$ such that $u_j\geq u$ and
	$$H_m(u_j)=\mathds{1}_{\{u>j\rho\}}H_m(u).$$
Note that, we have  $H_m(u_j)$ vanishes on all $m-$polar subsets of $\Omega$ for every $j$.
	Moreover, by Theorem 3.6 in \cite{HP17}, we have
	$$H_m(u_j)\leq \mathds{1}_{\{u>j\rho\}}H_m(\max(u,j\rho)\leq H_m(\max(u,j\rho).$$
	Moreover, since $u_j\geq u$ and $u\in\mathcal{N}_m(\Omega)$, we also have $u_j\in\mathcal{N}_m(\Omega).$ By Corollary 5.8 in $\cite{T19}$, it follows that $u_j\geq\max(u,j\rho)\geq j\rho$, and so that $u_j\in\mathcal{E}_m^0(\Omega)$.
	We observe that $H_m(u_{j+1})\geq H_m(u_j)$ for every $j$. Again, by Corollary 5.8 in $\cite{T19}$, we have $u_j\geq u_{j+1}$ for every $j$.\\
Since $F: [-\infty, 0] \times \Om \longrightarrow\mathbb [-\infty,0]$ be a continuous function that satisfies the conditions (A) and (B), we infer that $\delta_t$ is an increasing continuous function in $t$. Thus, we obtain $F(u,z)\leq F(0,z)\leq \delta_0\leq 0$ for all $z\in\Omega, t\leq 0.$ Obviously, $\delta_0>\delta_{-\infty}\geq-\infty.$\\
 Now, in the case $\delta_0<0$ then we have  that
	\begin{equation}\label{e4.2}
	\int_{\Omega}-\delta_0 H_m(u)\leq \int\limits_{\Omega} [-F( u(z),z)]H_m(u)  <\infty.
	\end{equation}
In the case $\delta_0=0.$ Using the argument as in the proof of Theorem 2.3 in \cite{Be11}, we can choose a non-decreasing function $\tilde{F}(t)$ such that $\tilde{F}^{'}=\tilde{F}^{''}=0$ on $[-1,0]$, $\tilde{F}$ is convex on $(-\infty,-1)$ and $\tilde{F}(t)\geq \delta_t.$ We have 
$$dd^c\tilde{F}(u_1)=\tilde{F}^{''}(u_1)du_1\wedge d^cu_1+\tilde{F}^{'}(u_1)dd^cu_1.$$ Thus, we obtain $\tilde{F}(u_1)\in SH_m^-(\Omega).$ Moreover, it follows from $u_1\in L^{\infty}(\Omega)$ that $\tilde{F}(u_1)\in SH_m^-(\Omega)\cap L^{\infty}(\Omega).$ Therefore, we have $\tilde{F}(u_1)\in\mathcal{E}_m(\Omega).$ We claim that 
\begin{equation}\label{e4.2b}
	\int_{\Omega}-F(u_1) H_m(u)  <\infty.
\end{equation}
	Indeed, we have
	\begin{align}\label{e4.4}
		\begin{split}
		\int\limits_{\Omega} [-\tilde{F}( u_1(z))]H_m(u)&\leq \int\limits_{\Omega} [-F( u_1(z),z)]H_m(u)\\
	&\leq \int\limits_{\Omega} [-F( u(z),z)]H_m(u) \\
	&\leq \int\limits_{\Omega} [-F( u(z),z)]H_m(u) <\infty.
	\end{split}
	\end{align}
	Now we set $v:=\lim\limits_{j\rightarrow\infty}u_j$. We have $v\geq u$ and $$H_m(v)=\lim\limits_{j\rightarrow\infty}H_m(u_j)=H_m(u).$$ Then, applying Theorem 2.10 in \cite{Gasmi2} coupling with inequalities \eqref{e4.2}, \eqref{e4.2b}  we have $u=v$, and so that $u_j\searrow u$. Finally,
	by monotone convergence theorem we get
	\begin{align*}
		\int\limits_\Omega[-F( u_j(z),z)]H_m(u_j)&=\int\limits_\Omega[-F( u_j(z),z)]\mathds{1}_{\{u>j\rho\}}H_m(u)\\
		&\rightarrow\int\limits_\Omega[-F( u(z),z)]H_m(u).
	\end{align*}
	The proof is completed.
\end{proof} 
As a consequence of the above results we have the following converse to Proposition \ref{theo4.4}.
\begin{corollary} \label{coro4.3}
Let $F: [-\infty, 0] \times \Om \longrightarrow\mathbb [-\infty,0]$ be a continuous function 
that satisfies the conditions (A) and (B). 
If $u\in\mathcal{N}_m(\Omega)$ and $$\int\limits_\Omega[-F (u(z),z)]H_m(u)<+\infty$$
then we have $u\in\mathcal{E}_{m,F}(\Omega).$
\end{corollary}
\n Next, we study the property of Hessian measure of functions of the classes $\E_{m,F}(\Om)$ when $F(-\infty,z)\equiv -\infty$ for all $z\in\Om.$  
\begin{proposition}\label{th4.6}
	Let $F: [-\infty, 0] \times \Om \longrightarrow\mathbb [-\infty,0]$ be an upper continuous function 
	that satisfies the conditions (A) and (B). If $F(-\infty,z)\equiv -\infty$ for all $z\in\Om$ then we have $\mathcal{E}_{m,F}(\Om)\subset\mathcal{N}_{m}^a(\Om).$
\end{proposition}
\begin{proof}
	Let $(u_j)$ be a sequence in $\mathcal{E}_m^0(\Omega)$ decreasing to $u\in \mathcal{E}_{m,F}(\Om)$ such that
	\begin{equation}\label{e4.5}
		M:=\sup\limits_j\int\limits_\Omega[-F( u_j(z),z)]H_m(u_j)<\infty.
		\end{equation}
	Fix $t<0.$ For all $j\leq k,$ since $\{u_j\}$ is a decreasing sequence, we have $u_j\geq u_k.$ Thus, we infer that $u_k<t$ on $\{u_j<t\}.$ Since $F$ is an increasing function in the first variable, we deduce that $F(u_k,z)\leq F(t,z)$  on $\{u_j<t\}.$ This is equivalent to $-F(u_k,z)\geq -F(t,z)$  on $\{u_j<t\}.$ Therefore, we obtaion 
	\begin{equation}\label{e4.6b}
		\int_{\{u_j<t\}}H_m(u_k)\leq \int_{\{u_j<t\}}\dfrac{-F(u_k,z)H_m(u_k)}{-F(t,z)}\leq   \int_{\Om} \dfrac{-F(u_k,z)H_m(u_k)}{-F(t,z)}.
	\end{equation}
Note that we have $F(t,z)\leq \delta_t<0.$ This implies that $-F(t,z)\geq -\delta_t>0.$ Thus, we have

 \begin{equation}\label{e4.6c}\dfrac{1}{-F(t,z)}\leq \dfrac{1}{-\delta_t}.
	\end{equation}
	Combining inequalities \eqref{e4.6b} and \eqref{e4.6c} we have

	\begin{equation}\label{e4.6d}
		\int_{\{u_j<t\}}H_m(u_k)\leq \int_{\Om} \dfrac{-F(u_k,z)H_m(u_k)}{-\delta_t}= \dfrac{1}{-\delta_t}\int_{\Om}-F(u_k,z)H_m(u_k).
	\end{equation}
	Combining inequalities \eqref{e4.5} and \eqref{e4.6d} we have
	\begin{equation}\label{e4.7}
		\int_{\{u_j<t\}}H_m(u_k)\leq \dfrac{1}{-\delta_t} M.
		\end{equation}
On the other hand, since $u_j$ is an upper semicontinuous function, we infer that $\{u_j<t\}$ is an open subset. This implies that $\ind_{\{u_j<t\}}$  is a lower semicontinuous function. Moreover, by Theorem \ref{theo4.1}, we have $u\in \mathcal{E}_{m,F}(\Om)\subset\mathcal{N}_m(\Om).$ Thus, by Theorem 3.8 in \cite{HP17}, we get $H_m(u_k)\to H_m(u)$ weakly as $k\to\infty.$  It follows from Lemma 1.9 in \cite{De93} that
\begin{equation}\label{i4.8}
	\liminf\limits_{k\to\infty}\int_{\{u_j<t\}}H_m(u_k)\geq \int_{\{u_j<t\}}H_m(u).
\end{equation}
Combining inequalities \eqref{e4.7} and \eqref{i4.8} we obtain
\begin{equation*}\label{e4.8}
	\int_{\{u_j<t\}}H_m(u)\leq \dfrac{1}{-\delta_t} M.
\end{equation*}
Letting $j\to\infty,$ by Lebesgue' monotone convergence Theorem we have
$$\int_{\{u<t\}}H_m(u)\leq \dfrac{1}{-\delta_t} M. $$

\n Note that the condition $F(-\infty,z)\equiv -\infty$ for all $z\in\Om$  is equivalent to $\delta_{-\infty}=-\infty.$ Now, letting $t\to-\infty$ we get
$\int_{\{u=-\infty\}}H_m(u)=0.$ This means that $u\in\mathcal{N}_m^a(\Om).$ The proof is complete.
\end{proof}

\subsection{ Complex $m$-Hessian equations in the class $\mathcal{E}_{m,F}(\Omega)$}
In \cite{HP17} the authors solved complex $m$-Hessian equation if it has subsolution in the class $\mathcal{E}_{m}(\Omega).$ In \cite{Gasmi}, A. El. Gasmi solved complex $m$-Hessian equation in the class $\mathcal{N}_{m}(f)$ if it has subsolution in the class $\mathcal{N}_{m}(\Omega).$ In \cite{HQ21}, L. M. Hai and V. V. Quan  solved $m-$Hessian type equation $H_{m}(u)=F(u,.)d\mu$ in the class $\mathcal{D}_{m}(\Omega)$ if there exists subsolution in this class. Next, in \cite{AAG20}, H. Amal, S. Asserda and A. El. Gasmi  solved the above mentioned equation  in the class $\mathcal{N}_{m}(f)$ if there exists subsolution in the class $\mathcal{N}_{m}(\Omega)$. Continuing the direction of the  above authors, in this Subsection, the authors will solve complex $m$-Hessian type equation $-F(u(z),z)H_{m}(u)=\mu$ in the class $\mathcal{E}_{m,F}(\Omega)$ in the case when measure $\mu$ vanishes on all of $m$-polar sets and in the case when  measure $\mu$ is arbitrary.\\
In this section, we let $F: [-\infty, 0] \times \Om \longrightarrow\mathbb [-\infty,0]$
be a {\it continuous} function such that $F: [-\infty, 0] \times \Om \longrightarrow\mathbb (-\infty,0]$ or $F(-\infty,z)\equiv -\infty$ for all $z\in\Om$ and $F$ satisfies the conditions (A) and (B) in Subsection 4.1.\\
As a typical function having the above properties we can take $F (t,z)=\psi(t) u(z)$ 
where $\psi: [-\infty, 0] \to (0,+\infty)$ is continuous decreasing and $u:\Om \to (-\infty, 0)$ is 
continuous.

\n 
First, we solve the complex $m$-Hessian equations in the class $\mathcal{E}_{m,F}(\Omega)$ when 
$\mu$ puts no mass on $m$-polar sets.
\begin{theorem}\label{th4.1}
 Let $\mu$ be a non-negative finite measure on $\Omega$ such that $\mu$ vanishes on all $m-$polar sets. Then there exists a unique function $u
\in \mathcal{N}^{a}_{m}(\Omega) \cap \mathcal{E}_{m,F}(\Omega)$ such that $H_{m,F}(u)=\mu.$
\end{theorem}
\begin{proof}
Using Theorem 3.1 in \cite{AAG20} in this case $f\equiv 0$, we can find a unique solution $u\in\mathcal{N}^{a}_{m}(\Omega)$ such that 
$$H_{m,F}(u)=-F(u(z),z)H_{m}(u)= \mu.$$
By the hypothesis $\mu(\Omega)<+\infty,$ we have $\int_{\Omega}H_{m,F}(u)<+\infty$. By Corollary \ref{coro4.3} and $u\in\mathcal{N}^{a}_{m}(\Omega)$, we conclude that $u\in\mathcal{E}_{m,F}(\Omega).$ The uniqueness of $u$ follows from Theorem \ref{comparison}.
The proof is complete.
\end{proof}
\n Next, we solve the complex $m$-Hessian equations in the class $\mathcal{E}_{m,F}(\Omega)$ when measure $\mu$ is dominated by certain weighted Hessian measure.
\begin{theorem}\label{theorem4.2}
Let $\mu$ be a non-negative finite measure on $\Omega$	such that there exists a function $w\in\mathcal E_{m,F}(\Omega)$ with 
$$\mu\leq -F(w(z),z)H_{m}(w).$$ 
Then there exists a function $u\in\mathcal{E}_{m,F}(\Omega)$ such that $u\geq w$ and $$-F(u(z),z)H_{m}(u)=\mu.$$
\end{theorem}
\begin{proof} We split the proof into two cases.\\
	{\bf Case 1.} $F: [-\infty, 0] \times \Om \longrightarrow\mathbb (-\infty,0].$
By the considering the equation $H_m(u)=-\frac{1}{F(u(z),z)}\mu$ and using the Main Theorem in \cite{AAG20} for $f=0,$ there exists $u\in\mathcal{N}_m(\Omega)$ such that $\mu=-F(u(z),z)H_m(u)$ with $u\geq w.$ By the hypothesis $w\in\mathcal{E}_{m,F}(\Omega)$ and Proposition \ref{theo4.4}	we have $\int_{\Omega}[-F(w(z),z)]H_m(w)<+\infty$. Thus, we have 
$$\int_{\Omega}[-F(u(z),z)]H_m(u) \leq \int_{\Omega}[-F(w(z),z)]H_m(w)<+\infty.$$
According to Corollary \ref{coro4.3} we get $u\in\mathcal{E}_{m,F}(\Omega).$\\
{\bf Case 2.} $F(-\infty,z)\equiv -\infty$ for all $z\in\Om.$\\
	By Proposition \ref{th4.6}, we have $w\in\E_{m,F}(\Om)\subset\mathcal{N}_m^a(\Om).$ It follows from $$\mu\leq -F(w,z)H_m(w)$$ that $\mu$ vanishes on all $m$-polar sets. Moreover, by Proposition \ref{theo4.4} we also have $$\mu(\Om)\leq\int_{\Om}-F(w,z)H_m(w)<+\infty.$$ Now, applying Theorem \ref{th4.1}, there exists a unique function $u\in\E_{m,F}(\Om)$ such that $\mu=-F(u,z)H_m(u).$ By Theorem \ref{comparison} we infer that $u\geq w.$
	The proof is complete.
\end{proof}

\section{Subextension in the class $\mathcal{E}_{m,F}$}

 Firstly, we need the following Lemma.
\begin{lemma}\label{keylemma}
	Let $F: [-\infty, 0] \times \Om \longrightarrow\mathbb (-\infty,0]$
	be a {\it continuous} function satisfying the conditions (A), (B) in Subsection 4.1 and let $\alpha$ be a finite Radon measure on $\Omega$ which puts no mass on $m-$polar sets of $\Omega$ such that $\text{\rm supp}\alpha\Subset\Omega$. Assuming that $v\in \mathcal F_{m}(\Omega)$ satifying $\text{\rm supp}H_{m}(v)\Subset\Omega$ and $H_{m}(v)$ is carried by a $m-$polar. We put
	$$u=(\sup \left\{\varphi: \varphi\in \mathcal U(\alpha, v)\right\})^*,$$
	where
	$$\mathcal U(\alpha, v)=\left\{\varphi\in SH^{-}_{m}(\Omega): \alpha\leq -F (\varphi(z),z)H_{m}(\varphi) \text{ and } \varphi\leq v \right\}.$$
	 Then we have $u\in\mathcal{F}_m(\Omega)$ and 
	$$-F(u(z),z)H_m(u) = \alpha -F (v(z),z) H_{m}(v).$$
\end{lemma}
\begin{proof}
		We split the proof into two steps as follows.
	
	\n 	
	{\em Step 1.}
	We prove that $u\in\mathcal{F}_m(\Omega)$ and
	\begin{equation}\label{eq5}
		-F (u(z),z)H_{m}(u)\geq\alpha -F(v(z),z)H_{m}(v).
	\end{equation}
	Indeed, by Theorem $\ref{th4.1}$ there exists a function $\phi\in\mathcal E_{m,F}(\Omega)\subset\mathcal{N}_m(\Omega)$ such that 
	$$ -F(\phi(z),z)H_{m}(\phi)=\alpha.$$
	\n
	Note that supp$H_{m}(\phi) \Subset \Omega$ so $\int_\O H_{m}(\phi)<+\infty$. Thus, by Theorem 4.9 in \cite{T19} we infer that $\phi\in\mathcal{F}_m(\Omega).$
	Moreover, it is easy to see that $(\phi+v)$ belongs to $\mathcal U(\alpha,v).$ So we have  $$\phi+v\leq u\leq v.$$ It follows that $u\in\mathcal{F}_m(\Omega).$ \\
	\n Next, using Choquet Lemma we may find a sequence $\{f_j\}\subset \mathcal U(\alpha,v)$ such that
	$$u=\Bigl(\sup_{j\in \mathbb N^*} f_j\Bigl)^*.$$
	Note that, by Lemma 3.1 in \cite{HQ21} we have $\max (\varphi,\psi)\in\mathcal U(\alpha, v)$, $\forall\ \varphi,\psi\in \mathcal U(\alpha, v)$. Thus, 
	$$\tilde u_j:=\max\{f_1,\ldots,f_j\} \in \mathcal U(\alpha,v)$$ and $\tilde u_j\nearrow u$ in $C_m$-capacity. By Corollary $\ref{nicecap}$ we also have $$-F (\tilde u_j(z),z)H_{m}(\tilde{u_{j}})\to -F (u(z),z)H_{m}(u)\,\,\text{weakly.} $$ It follows that $$-F (u(z),z)H_{m}(u)\geq\alpha.$$ Thus, we get  $u\in\mathcal U(\alpha,v)$.\\
	Now, since $\alpha$ puts no mass on $m-$polar sets, we have
	\begin{align}\label{e5.2}
		\begin{split}
			-F(u(z),z)H_{m}(u) &= -F(u(z),z) \ind_{\{u>-\infty\}} H_{m}(u) -F(u(z),z) \ind_{\{u=-\infty\}} H_{m}(u)\\
			&\geq  \ind_{\{u>-\infty\}} \al -F(u(z),z) \ind_{\{u=-\infty\}} H_{m}(u)\\
			&\geq \al -F(u(z),z) \ind_{\{u=-\infty\}} H_{m}(u).	
		\end{split}
	\end{align}
	\n On the other hand, since $H_{m}(v)$ is carried by an $m$-polar set then we have
	$$H_{m}(v)=\ind_{\{v>-\infty\}}H_{m}(v)+ \ind_{\{v=-\infty\}}H_{m}(v)=\ind_{\{v=-\infty\}}H_{m}(v).$$
	By Proposition 5.2 in \cite{HP17} we get
	\begin{equation}\label{e5.3}	 1_{\{u=-\infty\}} H_{m}(u)\geq 1_{\{v=-\infty\}} H_{m}(v) =H_{m}(v).
	\end{equation}
	Combining inequalities \eqref{e5.2} and \eqref{e5.3} we get
	$$-F(u(z),z)H_m(u)\geq \alpha -F(v(z),z) H_{m}(v).$$
	
	{\em Step 2.} We prove that
	$$-F(u(z),z)H_m(u) = \alpha -F (v(z),z) H_{m}(v).$$
	Choose a $m-$hyperconvex domain $G\Subset\O$ such that 
	$$\text{\rm supp}\alpha \cup\text{\rm supp}H_{m}(v)\Subset G\Subset \O.$$ Next, by $v\in \mathcal{F}_{m}(\O),$ we can choose 
	$\mathcal E_{m}^{0}(\Omega) \supset \{v_j\}\downarrow v $ such that 
	$$\sup_{j}\int_{\Omega}H_{m}(v_{j})<+\infty, \text{\rm supp}H_{m}(v_{j})\subset\overline{G}.$$ 
	Now since
	$$\int_\Omega \left[\alpha -F(v_j(z),z)H_{m}(v_{j})\right]\leq\int_{\overline{G}} \left[\alpha +\max_{\overline{G}}[-F(-\infty,z) ]H_{m}(v_{j})\right]<+\infty$$
	so by Theorem  $\ref{th4.1}$  there exists a unique function $w_j\in\mathcal E_{m,F}(\Omega)$ such that
	\begin{equation} \label{eq7}
		-F(w_j(z),z) H_{m}(w_{j}) =\alpha - F( v_j(z),z )H_{m}(v_{j}).
	\end{equation}
	On the other hand, we have
	$$ -F((\phi+v_j)(z),z)H_m(\phi)\geq -F(\phi(z),z)H_m(\phi)=\alpha$$
	and 
	$$ -F((\phi+v_j)(z),z)H_m(v_j)\geq -F(v_j(z),z)H_m(v_j).$$
	Thus, we have 
	
	\(
	\begin{aligned}[t]
		-F ((\phi+v_j)(z),z) H_m (\phi+v_j)&\geq\alpha -F(v_j(z),z)H_m(v_j)\\
		&= -F(w_j(z),z) H_{m}(w_j)\\
		&\geq -F (v_j(z),z) H_{m}(v_j).
	\end{aligned}
	\)
	
	\n By Theorem $\ref{comparison}$ we get $\phi+v_j\leq w_j\leq v_j.$ 
	It follows that $w_j\in\mathcal{F}_m(\Omega).$ Now we put 
	$$\mathcal V(\alpha, v_j)=\left\{\varphi\in \mathcal E_{m}(\Omega): \alpha\leq -F (\varphi,z)H_{m}(\varphi) \text{ and } \varphi\leq v_j \right\}$$
	and   $$u_j=\Bigl(\sup \left\{\varphi: \varphi\in \mathcal V(\alpha,v_j)\right\}\Bigl)^*,$$
	then we have $w_j\in\mathcal V(\alpha, v_j).$
	Obviously, $u_j$ is decreasing and $u_j\geq w_j$ for all $j\geq 1$. Note that, by $\mathcal{F}_m(\Omega)\ni w_j\leq u_j,$ we also have $u_j\in\mathcal{F}_m(\Omega).$ We will prove that $\{u_j\}$  decreases pointwise to $u$, as $j\to+\infty$. Indeed, assume that $u_j\searrow \tilde{u}$ as $j\to +\infty.$\\
	\n (i) We have $u\leq v\leq v_j.$ Thus, $u\in \mathcal{V}(\alpha,v_j).$ It follows that $u\leq u_j.$ Let $j\to+\infty$ we get that $u\leq \u.$\\
	\n (ii) On the other hand, we have $u_j\leq v_j.$ Let $j\to+\infty$ we obtain that $\u\leq v.$ Note that, since $u_j\searrow \u$ we get that $u_j$ converges to $\u$ in $Cap_m$-capacity. By Corollary $\ref{nicecap}$ we get that $-F(u_j(z),z)H_m(u_j) \to -F(\u(z),z)H_m(\u)\,\,\text{weakly} .$
	Therefore $-F(\u(z),z)H_m(\u)\geq \alpha.$ So we have $\u\in\mathcal{U}(\alpha,v)$ and we obtain $\u\leq u.$\\ 
	From (i) and (ii) we get $\u=u.$ That means	$\{u_j\}$  decreases pointwise to $u$, as $j\to+\infty$.\\
	Next,  we will prove that $u_j-w_j\to 0$ in $\tilde{C}_{m}$-capacity. Indeed, fix a strictly $m$-subharmonic function $h_0\in \mathcal E_m^{0}(\O)\cap \mathcal C^\infty(\O)$, such a function $h_0$ exists by \cite{AL}. For  $\varepsilon>0$ and $j_0\ge 1$ we have
	\begin{align*}
		& \limsup_{j\to\infty} \tilde{C}_{m} (\{u_j-w_j>\varepsilon\})
		\\&=\limsup_{j\to\infty}\Bigl( \sup\Bigl\{\int_{\{u_j-w_j>\varepsilon\}} dd^c h_0\wedge (dd^c h)^{m-1}\wedge\beta^{n-m}:
		\\&{~~~~~~~~~~~~~~~~~~~~~~~~~~~~~~~~~~~~~~~~~} h\in SH_m(\O), -1\leq h\leq 0\Bigl\}\Bigl)
		\\&\leq\limsup_{j\to\infty}\Bigl( \sup\Bigl\{\frac 1 { \varepsilon^m }\int_{\{u_j-w_j>\varepsilon\}} (u_j-w_j)^m  dd^c h_0\wedge (dd^c h)^{m-1}\wedge\beta^{n-m}:
		\\&{~~~~~~~~~~~~~~~~~~~~~~~~~~~~~~~~~~~~~~~~~} h\in SH_m(\O), -1\leq h\leq 0\Bigl\}\Bigl)
		\\&\leq\limsup_{j\to\infty}\Bigl( \sup\Bigl\{\frac 1 { \varepsilon^m }\int_{\O} (u_j-w_j)^m  dd^c h_0\wedge (dd^c h)^{m-1}\wedge\beta^{n-m}:
		\\&{~~~~~~~~~~~~~~~~~~~~~~~~~~~~~~~~~~~~~~~~~} h\in SH_m(\O), -1\leq h\leq 0\Bigl\}\Bigl)\\
		&\leq\limsup_{j\to\infty} \frac {m!} { \varepsilon^m } \int_{\O} -h_0 \Bigl[H_m(w_j)-H_m(u_j)\Bigl]
		\\&\leq\limsup_{j\to\infty} \frac {m!} { \varepsilon^m } \int_{\O} -h_0 H_m(w_j) +\frac {m!} { \varepsilon^m } \int_{\O} h_0 H_{m}(u)
		\\&\leq\limsup_{j\to\infty} \frac {m!} { \varepsilon^m } \int_{\O} -h_0 \frac {-F (w_j(z),z) H_m(w_j)} {-F (u_j(z),z)} +\frac {m!} { \varepsilon^m } \int_{\O} h_0 H_m(u)
		\\&=\limsup_{j\to\infty} \frac {m!} { \varepsilon^m } \int_{\O} -h_0 \frac {\alpha -F(v_j(z),z) H_m(v_j)} {-F (u_j(z),z)} +\frac {m!} { \varepsilon^m } \int_{\O} h_0 H_m(u)
		\\&\leq\limsup_{j\to\infty} \frac {m!} { \varepsilon^m } \int_{\O} -h_0 \frac {\alpha -F(v_j(z),z) H_m(v_j)} {-F (u_{j_0}(z),z)} +\frac {m!} { \varepsilon^m } \int_{\O} h_0 H_m(u),
	\end{align*}
	where the fourth estimate follows from Lemma 2.9 in \cite{Gasmi2} and the sixth estimate follows from supp$H_m(w_j)\Subset\Omega.$ Since $v_j\searrow v,$ we infer that $v_j\to v$ in $Cap_m$-capacity. Thus, it follows from Theorem \ref{cap} that $$-F(v_j(z),z)H_m(v_j)\to -F(v(z),z)H_m(v)\,\,\text{ weakly as}\,\, j\to+\infty.$$ Moreover, $F(u_{j_{0}}(z),z)$ is an upper semicontinuous function, so by Lemma 1.9 in \cite{De93} we have 
	$$\lim\limits_{j\to\infty}\frac { -F(v_j(z),z)H_m(v_j)} {-F(u_{j_0}(z),z)}\leq \frac {-F(v(z),z) H_m(v)} {-F (u_{j_0}(z),z)}.$$
	\n Note that $\text{\rm supp}\alpha \cup\text{\rm supp}H_m(v_j)\Subset\overline{G}$, we obtain
	\begin{align*}
		&\limsup_{j\to\infty}  \int_{\O} -h_0 \frac {\alpha -F(v_j(z),z) H_m(v_j)} {-F(u_{j_0}(z),z)} \\
		&= \limsup_{j\to\infty}  \int_{\overline{G}} -h_0 \frac {\alpha -F(v_j(z),z) H_m(v_j)} {-F (u_{j_0}(z),z)} \\
		&\leq  \int_{\overline{G}} -h_0 \frac {\alpha -F(v(z),z)  H_{m}(v)} {-F (u_{j_0}(z),z)} .\\
		&\leq  \int_{\O} -h_0 \frac {\alpha -F(v(z),z) H_{m}(v)} {-F(u_{j_0}(z),z)}. 
	\end{align*}
	Thus, we have
	\begin{align*}
		& \limsup_{j\to\infty} \tilde{C}_{m} (\{u_j-w_j>\varepsilon\})\\
		&\leq \frac {m!} { \varepsilon^m } \int_{\O} -h_0 \frac {\alpha -F(v(z),z) H_m(v)} {-F (u_{j_0}(z),z)} +\frac {m!} { \varepsilon^m } \int_{\O} h_0 H_m(u).
	\end{align*}	
	
	\n Letting $j_0\to\infty$ by the Lebegues monotone convergence Theorem and Step 1 we get
	\begin{align*}
		&\limsup_{j\to\infty} \tilde{C}_{m} (\{u_j-w_j>\varepsilon\})\\
		&\leq\limsup_{j_0\to\infty} \frac {m!} { \varepsilon^m } \int_{\O} -h_0 \frac {\alpha -F(v(z),z) H_m(v)} {-F (u_{j_0}(z),z)} +\frac {m!} { \varepsilon^m } \int_{\O} h_0 H_m(u)
		\\&= \frac {m!} { \varepsilon^m } \int_{\O} -h_0 [\frac {\alpha -F(v(z),z) H_m(v)} {-F(u(z),z)} - H_m(u)]\leq 0.
	\end{align*}
	Thus we have proved that $u_j-w_j \to 0$ in $\tilde{C}_{m}$-capacity.
	Since $u_j \to u$ in $Cap_m$-capacity, we infer that $w_j\to u$ in $\tilde{C}_{m}$-capacity. Hence by Proposition \ref{weakcap} we get
	$$\liminf\limits_{j\to\infty} -F(w_j(z),z) H_m(w_j)\geq -F(u(z),z) H_m(u).$$
	This yields that $$\alpha -F(v(z),z)H_m(v)\geq -F(u(z),z) H_m(u).$$
	Combining this with \eqref{eq5}, we finally reach
	$$-F(u,z) H_m(u)=\alpha -F(v,z) H_m(v).$$ The proof is complete.
\end{proof}

\n Now, we will prove the main result in this Section.
\begin{theorem}
	Let $\Omega\Subset \widetilde{\Omega}$ be bounded $m$-hyperconvex domains in $\mathbb{C}^{n}.$  Assume that $F: [-\infty, 0] \times \tilde{\Om} \longrightarrow\mathbb [-\infty,0]$
	is a {\it continuous} function such that $F: [-\infty, 0] \times \tilde{\Om} \longrightarrow\mathbb (-\infty,0]$ or $F(-\infty,z)\equiv -\infty$ for all $z\in\tilde{\Om}$ and $F$ satisfies the following conditions \\
	\n (C) For $z\in\tilde{\Omega}, t\mapsto F(t,z)$ is increasing in t.\\
	\n (D) For all $t<0$ given, there exists a constant $\delta_t<0$ (depend only $t$) such that $\sup\limits_{z\in\tilde{\Omega}}F(t,z) =\delta_t<0.$\\
	Then for every $u\in\mathcal{E}_{m,F}(\Omega)$ that satisfies $\int_{\Omega}H_m(u)<+\infty,$ there exists $\tilde{u}\in\mathcal{E}_{m,F}(\widetilde{\Omega})$ such that $\tilde{u}\leq u$ on $\Omega$ and 
	$$-F(\tilde{u},z)H_{m}(\tilde{u})=\ind_{\Omega}[-F(u,z]H_m(u)\,\,\text{on}\,\, \widetilde{\Omega}.$$ 
\end{theorem}
\begin{proof}
	We split the proof into two cases.\\
	{\bf Case 1.} $F: [-\infty, 0] \times \tilde{\Omega} \longrightarrow\mathbb (-\infty,0].$
	 We have $$-\infty<F(-\infty,z)\leq\delta_{-\infty}=\sup\limits_{z\in\tilde{\Omega}}F(-\infty,z)<0.$$ This implies that 
	\begin{equation}\label{equation5.7}0<\dfrac{1}{-F(-\infty,z)}\leq \dfrac{1}{-\delta_{-\infty}}<\infty.
	\end{equation} 
	By considering $-F (t,z)/F(-\infty,z)$ instead of $F$ and using the inequality \eqref{equation5.7},
	we can assume that $$F(-\infty,z) =-1, \forall z\in\tilde{\Omega}.$$
	It follows from $u\in\mathcal{E}_{m,F}(\Omega)$ and Theorem $\ref{theo4.1}$ that we have $u\in\mathcal{N}_{m}(\Omega).$ By hypothesis $\int_{\Omega}H_m(u)<+\infty$ and Theorem 4.9 in \cite{T19} we get $u\in\mathcal{F}_{m}(\Omega).$ Obviously we also have $$\ind_{\{u=-\infty\}}H_m(u)\leq H_m(u).$$ According to the main Theorem in $\cite{Gasmi}$ implies that there exists a function $v\in \mathcal{N}_{m}(\Omega)$ such that 
	\begin{equation}\label{eqe5.7}H_m(v)=\ind_{\{u=-\infty\}}H_m(u)
	\end{equation}
	and $v\geq u.$ Since $u\in\mathcal{F}_m(\Omega,$ we get $v\in\mathcal{F}_m(\Omega).$\\
	Let $\tilde{v}$ be the maximal subextension of $v$ on $\widetilde{\Omega}$  defined by
	$$\tilde{v}=\sup\{\varphi\in SH_m^{-}(\widetilde{\Omega}): \varphi\leq v \, \text{on}\,  \Omega\}.$$
	By repeating the arguments as in Theorem 3.3 of \cite{PDmmj} (when $\chi\equiv -1$), 
	we infer that $\tilde{v}\in\mathcal{F}_m(\widetilde{\Omega})$ and
	\begin{equation}\label{e5.7} \int_{\widetilde{\Omega}}H_m(\tilde{v})\leq \int_{\Omega}H_m(v)
		\leq \int_{\widetilde{\Omega}}\ind_{\Omega}H_m(v).
	\end{equation}
	Note that $H_m(v)$ is carried by an $m$-polar set. Moreover, we have $\tilde{v}\leq v$ on $\Omega$, so by Proposition 5.2 in \cite{HP17} we obtain
	$$H_m(v)= \ind_{\{v=-\infty\}}H_m(v)\leq \ind_{\{\tilde{v}=-\infty\}}H_m(\tilde{v}) \leq H_m(\tilde{v})\,\,\text{on}\,\, \Omega.$$ Thus, we have $H_m(\tilde{v})\geq \ind_{\Omega}H_m(v) $ on $\widetilde{\Omega}.$ Coupling with inequality \eqref{e5.7} we obtain that 
	\begin{equation}\label{eqe5.8}H_m(\tilde{v})=\ind_{\Omega}H_m(v) \,\,\text{ on}\,\, \widetilde{\Omega}.
	\end{equation} That means we have $H_m(\tilde{v})$ is carried by a $m$-polar set and supp$H_m(\tilde{v})\Subset \widetilde{\Omega}.$\\
	Now, we set 
	$$\tilde{u}= (\sup\{\varphi:\varphi\in \mathcal{U}( \tilde{\alpha},\tilde{v}) \})^{*},$$ where $\tilde{\alpha} =\ind_{\Omega\cap\{u>-\infty\}}[-F(u(z),z)]H_m(u)$ is a finite Radon measure on $\widetilde{\Omega}$ which put no mass on $m$-polar sets on $\widetilde{\Omega}$ such that supp$\tilde{\alpha}\Subset\widetilde{\Omega}$  and $$\mathcal{U}( \tilde{\alpha},\tilde{v})= \{\varphi \in\mathcal{E}_m(\widetilde{\Omega}): \tilde{\alpha}\leq -F(\varphi(z),z)H_m(\varphi)\,\text{and}\,\varphi\leq\tilde{v} \}.$$  
	By Lemma \ref{keylemma} we get that $\tilde{u}\in\mathcal{F}_m(\widetilde{\Omega})$ and coupling with equalities \eqref{eqe5.7} and \eqref{eqe5.8} we also have
	\begin{align*}&-F(\tilde{u}(z),z)H_m(\tilde{u})\\
		&= \ind_{\Omega\cap\{u>-\infty\}}[-F(u(z),z)]H_m(u) + [-F(\tilde{v}(z),z)]H_m(\tilde{v})\\
		&=\ind_{\Omega\{u>-\infty\}}[-F(u(z),z)]H_m(u) +[-F(\tilde{v}(z),z)] \ind_{\{\tilde{v}=-\infty\}}H_m(\tilde{v})\\
		&=\ind_{\Omega\{u>-\infty\}}[-F(u(z),z)]H_m(u) +[-F(\tilde{v}(z),z)] \ind_{\{\tilde{v}=-\infty\}}\ind_{\Omega}H_m(v)\\
		&=\ind_{\Omega\{u>-\infty\}}[-F(u(z),z)]H_m(u) +[-F(-\infty,z)] \ind_{\{\tilde{v}=-\infty\}}\ind_{\Omega}H_m(v)\\
		&=\ind_{\Omega\{u>-\infty\}}[-F(u(z),z)]H_m(u) +[-F(-\infty,z)] \ind_{\{\tilde{v}=-\infty\}}\ind_{\Omega}\ind_{\{v=-\infty\}}H_m(v)\\
		&=\ind_{\Omega\{u>-\infty\}}[-F(u(z),z)]H_m(u) +[-F(-\infty,z)] \ind_{\Omega}\ind_{\{v=-\infty\}}H_m(v)\\
		&=\ind_{\Omega\{u>-\infty\}}[-F(u(z),z)]H_m(u) +[-F(-\infty,z)] \ind_{\Omega}H_m(v)\\
		&= \ind_{\Omega\{u>-\infty\}}[-F(u(z),z)]H_m(u) +[-F(-\infty,z)] \ind_{\Omega}\ind_{\{u=-\infty\}}H_m(u)\\
		&= \ind_{\Omega\cap\{u>-\infty\}}[-F(u(z),z)]H_m(u)+\ind_{\Omega\cap\{u=-\infty\}}[-F(-\infty,z)]H_m(u)\\
		&= \ind_{\Omega\cap\{u>-\infty\}}[-F(u(z),z)]H_m(u)+\ind_{\Omega\cap\{u=-\infty\}}[-F(u(z),z)]H_m(u)\\
		&=\ind_{\Omega} [-F(u(z),z)]H_m(u),
	\end{align*}
	where the fifth equality follows from the fact that $H_m(v)$ is carried by an $m-$polar set and the sixth equality follows from the fact that $\tilde{v}\leq v $ on $\Omega.$
	Now, by Proposition \ref{theo4.4} we have $$\int_{\widetilde{\Omega}}-F(\tilde{u}(z),z)H_m(\tilde{u})=\int_{\Omega}[-F(u(z),z)]H_m(u)<+\infty.$$ It follows from Corollary \ref{coro4.3} that $\tilde{u}\in\mathcal{E}_{m,F}(\widetilde{\Omega}).$\\
	Thus, it remains to prove that $\tilde{u}\leq u$ on $\Omega.$
	Put $$\psi= (\sup\{\varphi:\varphi\in \mathcal{U}(\alpha ,v)\})^{*},$$ where $\alpha = \ind_{\{u>-\infty\}}[-F(u(z),z)]H_m(u)$ and 
	$$\mathcal{U}( \alpha,v)= \{\varphi \in\mathcal{E}_m(\Omega): \alpha\leq -F(\varphi(z),z)H_m(\varphi)\,\text{and}\,\varphi\leq v \}.$$  
	Note that $H_m(v)$ is carried by an $m$-polar set and we do not have supp$\alpha\Subset\Omega.$
	Thus, by repeating the argument as in Step 1 in the proof of Lemma \ref{keylemma}  we only have $\psi\in\mathcal{N}_{m}(\Omega)$ and  
	\begin{align*}&-F(\psi(z),z)H_m(\psi)\\
		&\geq \ind_{\{u>-\infty\}}[-F(u(z),z)]H_m(u) + [-F(v(z),z)]H_m(v)\\
		&=\ind_{\{u>-\infty\}}[-F(u(z),z)]H_m(u) + [-F(v(z),z)]\ind_{\{v=-\infty\}}H_m(v)\\
		&=\ind_{\{u>-\infty\}}[-F(u(z),z)]H_m(u) + [-F(-\infty,z)]\ind_{\{v=-\infty\}}H_m(v)\\
		&=\ind_{\{u>-\infty\}}[-F(u(z),z)]H_m(u) + [-F(-\infty,z)]H_m(v)\\
		&= \ind_{\{u>-\infty\}}[-F(u(z),z)]H_m(u) +[-F(-\infty,z)]\ind_{\{u=-\infty\}}H_m(u)\\
		&\geq \ind_{\{u>-\infty\}}[-F(u(z),z)]H_m(u)+\ind_{\{u=-\infty\}}[-F(u(z),z)]H_m(u)\\
		&\geq -F(u(z),z)H_m(u).
	\end{align*}
	Obviously, we have $u$ belongs to the class of functions which are in the definition of $\psi.$ Thus, $u\leq \psi.$ So we have
	$$-F(\psi(z),z)H_m(\psi)\geq -F(u(z),z)H_m(u)\geq -F(\psi(z),z)H_m(u). $$
	This implies that $H_m(\psi)\geq H_m(u)$ on $\Omega.$\\
	According to Lemma 2.14 in $\cite{Gasmi}$ and $u\leq \psi$ we obtain that
	$$+\infty>\int_{\Omega}H_m(u)\geq \int_{\Omega}H_m(\psi).$$
	So we have $H_m(\psi)=H_m(u).$ Theorem 2.10 in $\cite{Gasmi2}$ implies that $\psi=u.$ Now by the definition of $\tilde{u}$ and $u$ we get $\tilde{u}\leq u$ on $\Omega.$\\
	{\bf Case 2.} $F(-\infty,z)\equiv -\infty$ for all $z\in\tilde{\Omega}.$ By Proposition \ref{th4.6}, we have $u\in\E_{m,F}(\Om)\subset\mathcal{N}_m^a(\Om).$ Therefore, measure $\ind_{\Om}[-F(u,z)]H_m(u)$ vanishes on all $m$-polar subsets of $\tilde{\Om}$. Moreover, by Proposition \ref{theo4.4} we also have $$\int_{\tilde{\Om}}\ind_{\Om}[-F(u,z)]H_m(u)\leq\int_{\Om}-F(u,z)H_m(u)<+\infty.$$ By Theorem \ref{th4.1}, there exists a unique function $\tilde{u}\in\E_{m,F}(\tilde{\Om})$ such that $$-F(\tilde{u},z)H_m(\tilde{u})=\ind_{\Om}[-F(u,z)]H_m(u).$$ Now, it remains to prove that $\tilde{u}\leq u$ on $\Om.$ Indeed, we put $w=\max(u,\tilde{u})$ on $\Om.$ By Lemma 3.1 in \cite{HQ21} with measure $\mu=-F(u,z)H_m(z)$ puts no mass on all $m-$polar sets, we have 
	$$-F(w,z)H_m(w)\geq -F(u,z)H_m(u) \,\,\text{on}\,\,\Om.$$ By Theorem \ref{comparison} we infer that $u\geq w$ on $\Om.$ Thus, we obtain $\tilde{u}\leq u$ on $\Om.$ The proof is complete.
\end{proof} 

\section*{Declarations}
\subsection*{Ethical Approval}
This declaration is not applicable.
\subsection*{Competing interests}
The authors have no conflicts of interest to declare that are relevant to the content of this article.
\subsection*{Authors' contributions }
Hoang Thieu Anh, Nguyen Van Phu and Nguyen Quang Dieu and  together studied  the manuscript.
\subsection*{Availability of data and materials}
This declaration is not applicable.

\end{document}